\documentclass[20pt]{article}
\usepackage[margin=1in]{geometry}
\usepackage{dsfont}
\usepackage{amsmath, amsthm, amssymb, mathrsfs}
\usepackage[abbrev,non-sorted-cites]{amsrefs}
\usepackage[all]{xy}
\usepackage{graphicx}
\usepackage{extarrows}
\usepackage{hyperref}
\usepackage{xcolor}
\usepackage{tikz-cd}
\newtheorem{thm}{Theorem}

\newtheorem{lem}{Lemma}

\begin{document}
\begin{center}\begin{Huge}
Sums with Stern-Brocot sequences and  Minkowski question mark function 
\end{Huge}\end{center}
\begin{center}\begin{large}Liu Haomin,\quad Lü Jiadong,\quad Xie Yonghao\end{large}\end{center}

\vskip+1cm

\begin{large}\begin{bf}Abstract\end{bf}\end{large}\,\,
Motivated by a famous Franel's theorem dealing with distributions of Farey series, we consider an analogous problem related to Stern Brocot sequences. We prove that the remainder
$$
R_n=\sum_{j=1}^{2^n}\left(\xi_{j,n}-\frac{j}{2^n}\right)^2-2^n\int_0^1(?(x)-x))^2\text{d}x
$$
tends to $0$ when $n$ tends to infinity, where $\xi_{j,n}$ are elements of the Stern-Brocot sequence and $?(x)$ denotes Minkowski Question-Mark Function. 
We present some extended results and give a correct proof of a theorem on the Fourier-Stieltjes coefficient of the inverse function of $?(x)$.

\vskip+2cm

\section{Distribution of rational numbers: objects and definitions}

In this introductory section we discuss the   definitions of Stern-Brocot sequences and Minkowski function $?(x)$, and formulate
several recent and classical results.
In particular,  we recall
the famous Franel's theorem about the distribution of rational numbers from $[0,1]$ with bounded denominators.

\subsection{Stern-Brocot sequences}

The inductive  definition of Stern-Brocot sequences  $ F_n $, $ n=0, 1, 2, \dots $ is as follows. 
For $n=0$  define
  $$ F_0=\{0, 1\}=\left\{\frac{0}{1}, \frac{1}{1}\right\} .$$
Suppose that  for $ n \ge 0$ the sequence $ F_n $ is written in the increasing order
 $$0=\xi_{0, n}<\xi_{1, n}< \dots <\xi_{{ 2^{n}}, n}=1, \,\,\,,
\,\,\,\, \xi_{j, n} =\frac{p_{j, n}}{q_{j, n}},\,\,\, \gcd( p_{j, n},q_{j, n}) = 1.$$
Then the sequence  $F_{n+1}$ is defined as   $$ F_{n+1} = F_n \cup
W_{n+1} $$ where 
$$
 W_{n+1}=\left\{
\frac{p_{j, n}+ p_{j+1, n} }{q_{j, n}+ q_{j+1, n}},\,\,\,
 j=0, \dots , { 2^{n}}-1\right\}.
   $$  
{ Two fractions
$$
\frac{p_{j-1, n}+ p_{j, n} }{q_{j-1, n}+ q_{j, n}},\,\,\,\,\,
\frac{p_{j, n}+ p_{j+1, n} }{q_{j, n}+ q_{j+1, n}} \in W_{n+1}
$$  are called the  descendant of
the fraction 
 $\frac{p_{j,n}}{q_{j,n}} \in F_n$ , $ 1\le j \le {2^{n}}-1$. }
It is also worth mentioning that (see the proof of  Lemma 2  from \cite{g}) 
\begin{equation}\label{wn+1}
W_{n+1}=\left\{\frac{q}{2q-p},\frac{q-p}{2q-p}:\frac{p}{q}\in W_n\right\},\,\,\, n\ge 1.
\end{equation}
Note that for the number of elements in $F_n$ one has
{
$
|F_n|   =2^{n}+1
$,
and we have symmetry
\begin{equation}\label{sYM}
\xi_{j,n}+ \xi_{2^n-j,n} =1,\,\,\,\,\ 0\le j \le 2^n.
\end{equation}
}
First five sequences $F_0, F_1, F_2, F_3, F_4$ are visualised 
in Figure 1.

\vskip+0.3cm
\begin{tikzcd}[row sep=0.7em, column sep=0.8em]
\frac{0}{1} & & & & & & & & & & & & & & & & \frac{1}{1} & & F_0\\
\frac{0}{1} & & & & & & & & \frac{1}{2} & & & & & & & & \frac{1}{1} & & F_1\\
\frac{0}{1} & & & &\frac{1}{3} & & & &\frac{1}{2} & & & &\frac{2}{3} & & & & \frac{1}{1} & & F_2\\
\frac{0}{1} & & \frac{1}{4} & & \frac{1}{3} & & \frac{2}{5} & & \frac{1}{2} & & \frac{3}{5} & & \frac{2}{3} & & \frac{3}{4} & & \frac{1}{1} & & F_3\\
\frac{0}{1} & \frac{1}{5} & \frac{1}{4} & \frac{2}{7} & \frac{1}{3} & \frac{3}{8} & \frac{2}{5} & \frac{3}{7}& \frac{1}{2} & \frac{4}{7} & \frac{3}{5} & \frac{5}{8} & \frac{2}{3} & \frac{5}{7} & \frac{3}{4} & \frac{4}{5} & \frac{1}{1} & & F_4\\
& & & & & & & & \begin{bf}Figure\,1\end{bf} & & & & & & & &
\end{tikzcd}

\vskip+0.3cm
For every rational number $\xi\in [0,1]\cap \mathbb{Q}$ there exists minimal $n$ and unique $j$  from the range $0\le j \le 2^n$ such that $ 
\xi = \xi_{j,n}$. 
By $ \xi =[0;a_1,a_2,...,a_t] ,\,\,\, a_j \in \mathbb{Z}_+, a_t \ge 2$ we denote the unique  decomposition of rational $\xi \in [0,1]\cap \mathbb{Q}$ as a finite ordinary continued fraction.

The set $W_n$   for  $n\ge 1$ can be  characterised in terms of sums of partial quotients of its elements as
\begin{equation}\label{wn}
W_{n}=\left\{ \xi =[0;a_1,a_2,...,a_t] \in [0,1]\cap \mathbb{Q}: \,\,\,
a_1+ a_2+ ...+a_t = n+1
 \right\}.
\end{equation}
Each rational number $\xi $ can be uniquely written in the form $ \xi = \frac{p}{q}, \,(p,q) =1.$
If $
\xi = \frac{p}{q} = \frac{p_{j, n}+ p_{j+1, n} }{q_{j, n}+ q_{j+1, n}}
,$
then 
$q\ge \max( q_{j, n}, q_{j+1, n})$ and 
\begin{equation}\label{kss0}
S_n:=
\sum_{\frac{p}{q} \in W_{n}} \frac{1}{q^2}  
\le \sum_{j=0}^{2^n} \frac{1}{ q_{j, n} q_{j+1, n}}
=
\sum_{j=0}^{2^n} \left(\frac{p_{j+1,n}}{ q_{j+1, n}}-  \frac{p_{j,n}}{ q_{j, n}}\right)
 =1.
\end{equation}
Consider the Lebesgue measure
$\lambda (\mathcal{W}_n)$ of the set of real numbers 
$$
\mathcal{W}_n = 
\left\{ \xi = [0;a_1,a_2,...,a_\nu,...]  \in [0,1] : \,\,\,\exists\, t\,\,
\text{such that}\,\,
a_1+ a_2+ ...+a_t = n+1
 \right\}.
 $$
M. Kesseböhmer and B.O. Stratmann \cite{ks}  (see also B. Heersink \cite{hr}) proved that  for the value $S_n$ defined in 
(\ref{kss0}) we have
\begin{equation}\label{kss}
S_n  \asymp
\lambda (\mathcal{W}_n) \sim 
\frac{1}{\log_2 n},\,\,\,
n \to\infty.
\end{equation}

\subsection{Minkowski question mark function}

The Minkowski question mark function $?(x)$  can be  defined  as the limit distribution function for the Stern-Brocot sequences by the formula
$$
?(x)  =\lim_{n\to \infty} \frac{| F_n \cap [0,x)|}{|F_n|}=\lim_{n\to \infty} \frac{| F_n \cap [0,x)|}{2^n +1}.
$$
It is well-known that $?(x)$ 
  is continuous strictly increasing
  singular function. In terms of continued fraction expansion of $ x \in [0,1]$ the value $?(x)$ may be written as 
  \begin{equation}
  \label{eqn-1}
  ?(x)=2\sum_{n=1}^\infty(-1)^{n-1}2^{-\sum_{k=1}^na_k},
  \end{equation}
  where $a_k$ denotes the $k$-th partial quotient of $x$ and the sum is finite in the case when $ x\in \mathbb{Q}$.
  Question mark function $?(x)$ satisfies identities 
 {\begin{equation}
 \label{eqn-2}
 ?(1-x)=1-?(x),\,\,\,\,\,\,\, ?\left(\frac{x}{1+x}\right)=\frac{?(x)}{2}.
 \end{equation}}
 The function $m(x)$ inverse to $?(x)$ is also a  monotone continuous singular function.
 
 All the  mentioned above basic facts about  the Minkowski  function  are discussed for example in papers \cite{a1,a2}.

\vskip+0.3cm
In 1943 Salem \cite{s}
 proposed the following problem concerning the Fourier-Stieltjes coefficients of $?(x)$:
  is it true that 
  $$\widehat{?}(n)=\int_0^1e(nx)\text{d}?(x)\to 0,\,\,n\to\infty\,\,\,\,\,?$$
  (We use standard notation
  $e(z) = e^{2\pi i z}$.)
  An affirmative answer was given by T. Jordan and T. Sahlsten in 2016 in \cite{js}.
  They proved this long standing conjecture and showed that  indeed $\widehat{?}(n)\to0,\,\,n\to\infty.$

\vskip+0.3cm
In 2014 E.P. Golubeva \cite{g} considered an analogue of Salem's problem for  Fourier-Stieltjes coefficients
$$
\widehat{m}(n)=
\int_0^1e(nx)\text{d}m(x)
$$  of the inverse  function $m(x)$.  The situation with the inverse function is quite different. Golubeva proved \cite{g} that
 $\widehat{m}(n)$ do  not tend to zero as $n$ tends to infinity.
Paper \cite{g} is very important for our consideration, because the constructions from the present paper rely essentially on the approach from \cite{g}.

\subsection{Farey sequences and Franel's theorem}
Here we consider Farey sequences (or Farey series) ${\cal F}_Q$ which consist of all rational numbers $p/q\in [0,1], (p,q) = 1$
with denominators $\le Q$. Suppose that ${\cal F}_Q$ form an increasing sequence
$$
1= r_{0,Q}<r_{1,Q}<...< r_{j,Q}<r_{j+1,Q} <...< r_{\Phi (Q),Q}= 1, \,\,\,\,
\Phi (Q) 
=\sum_{q\le Q}\varphi (q)
$$
(here $\varphi (\cdot )$ is the Euler totient function).
It is well known that 
\begin{equation}\label{a0}
\lim_{Q\to \infty} \frac{|{\cal F}_Q \cap[0,x)|}{| {\cal F}_Q| } =
\lim_{Q\to \infty} \frac{|{\cal F}_Q \cap[0,x)|}{\Phi (Q)+1} =x.
\end{equation}
The famous 
Franel's theorem (see \cite{f,l}) states that the asymptotic formula
$$
 \sum_{j=1}^{\Phi (Q)} \left(r_{j,Q} - \frac{j}{\Phi (Q)}\right)^2 = O_\varepsilon (Q^{-1+\varepsilon}) 
, \,\,\, Q\to \infty.
$$
for all positive $\varepsilon$ 
is equivalent to Riemann Hypothesis.
In fact the well-known asymptotic equality 
for M\"{o}bius function
$$
\sum_{n\le Q} \mu (n) = o(Q), \,\,\,
Q\to \infty
$$
(which is equivalent to Prime Number Theorem)
leads to
\begin{equation}\label{a1}
\sum_{j=1}^{\Phi (Q)} \left(r_{j,Q} - \frac{j}{\Phi (Q)}\right)^2 = o(1) 
, \,\,\, Q\to \infty.
\end{equation}
All the details one can find in a wonderful book by E. Landau (see  \cite{landau}, Ch. 13).

\section{Main results}

We divide our main results into three subsections.

\subsection{Distribution of Stern-Brocot sequences}

 Applying  Koksma's  inequality (see \cite{n}, Ch. 2, \S 5)  Moshchevitin \cite{m} showed that 
\begin{equation}\label{fraa}
 \sum_{j=1}^{2^n} \left(\xi_{j,n} - \frac{j}{2^n}\right)^2 =
2^n \int_0^1 ( ?(x) - x)^2 {\rm d}x +R_n,\,\,\,
|R_n|\le 4,\,\,\,  n=1,2,3,....\, . 
\end{equation}
Here the sum with elements of  Stern-Brocot sequence  in  left hand side  is similar to the expression from  the left hand side  of  (\ref{a1}) with Farey fractions.
Meanwhile in the right hand side of  (\ref{a1}) there is no main term, as  the integral with the distribution function analogous to 
$\int_0^1 ( ?(x) - x)^2 {\rm d}x $ is equal to zero, because of  (\ref{a0}). In the present paper we show that the remainder $R_n$ in the right hand side of (\ref{fraa}) tends to zero and prove the following result.

\begin{thm}\label{thm1}
When $n\to \infty$ we have
$$
\sum_{j=1}^{2^n}\left(\xi_{j,n}-\frac{j}{2^n}\right)^2=
2^n\int_0^1(?(x)-x)^2\text{\rm d}x +O\left(n^{-\frac{3}{2}}(\log n)^{-\frac{1}{2}}\right).
$$
\end{thm}

 \subsection{A general result and its corollaries}

We have the following general result.

\begin{thm} \label{thmBV}
 Let $j(x)$ be 1-periodic odd piecewise continuous function
 with 
 $$
 \sup_{x\in [0,1]} |j(x)| <\infty.
 $$
Let    $F(x)$   be a continuous function which  has bounded variation on $\left[\frac{1}{2},1\right]$.
Then
\begin{equation}\label{oo}
\int_{\frac{1}{2}}^1j(2^n?(x))F(x)\text{d}x= O(n^{-1}), \,\,\,n \to\infty.
\end{equation}
\end{thm}

From Theorem \ref{thmBV} we deduce several corollaries. First of all,
It should be noticed that 
 Theorem \ref{thmBV} leads to  the asypmtotic equality 
\begin{equation}\label{simpler}
 \sum_{j=1}^{2^n}\left(\xi_{j,n}-\frac{j}{2^n}\right)^2= 2^n\int_0^1(?(x)-x)^2\text{d}x+O(n^{-1}), \,\,\,n \to\infty.
\end{equation}
This statement  is weaker  than the result of our Theorem \ref{thm1}, however its proof is essentially simpler.
It is given in Section \ref{simp}.
 \vskip+0.3cm
 Then we deduce another general corollary.

 \vskip+0.3cm
 \noindent
{\bf Corollary 1.}  {\it 
Let $F(x)$ be continuous function with bounded variation on the segment
  $\left[\frac{1}{2},1\right]$.
  Then
$$ \sum_{k=2^{n-1}}^{2^n-1}(-1)^k\int_{\xi_{k,n}}^{\xi_{k+1,n}}F(x)\text{\rm d}x=
O(n^{-1}), \,\,\,n \to\infty.$$
}
 
  \vskip+0.3cm
 To illustrate the result of Corollary 1 we consider two further examples.
 {
 We should note that from symmetry property (\ref{sYM}) it follows that 
 $$
 \sum_{j=0}^{2^n} \xi_{j,n} = \frac{1}{2} + 2^{n-1},\,\,\,\,\,
 \sum_{j=0}^{2^n}(-1)^j \xi_{j,n} =\frac{1}{2} ,\,\,\,\,\, n \ge 1.
 $$
 Our corollaries deal with asymptotic formulas for  sums of higher powers of $\xi_{j,n}$.}

 \vskip+0.3cm
 \noindent
{\bf Corollary 2.} {\it For every nonnegative integer $m$ one has}
$$
 \sum_{k=2^{n-1}}^{2^n}(-1)^k(\xi_{k,n})^m=\frac{1}{2}+\frac{1}{2^{m+1}}
 + O(n^{-1}), \,\,\,n \to\infty.
$$

 \vskip+0.3cm
  \noindent
{\bf Corollary 3.} {\it For every positive integer $m$ one has }
$$
\sum_{k=0}^{2^n}(-1)^k(\xi_{k,n})^m=\frac{1}{2}
+O(n^{-1}), \,\,\,n \to\infty.
$$

  \vskip+0.3cm
  Corollaries 1, 2, 3 are proven in Subsection \ref{bbb}.

\subsection{Fourier coefficients for the inverse function}

This section is devoted to the behaviour of the  Fourier-Stieltjes coefficients 
$\widehat{m}(n)$
of the inverse function $m(x)$. In the introduction we  refereed to Golubeva's result \cite{g} which states that 
$$
\limsup_{n\to \infty}
\widehat{m}(n) > 0.
$$
In \cite{gg} Gorbatyuk claimed a stronger result that 
\begin{equation}\label{eqGoT}
\lim_{n\to\infty}\widehat{m}(2^n)=1.
\end{equation}
However, Gorbatyuk's proof contained a mistake. Gorbatyuk claims the  equality
\begin{equation}\label{eqGo}
T^n ( f ) = L^n ( f ),\,\,\,\ n =1,2,3,...
\end{equation}
for certain operators $T,L$\footnote{here we use original  notation     $T$ from  \cite{gg}  for a  certain operator 
which differs from our operator defined in (\ref{ope}).}
  on functions $f:[0,1]\to [ 0,1]$ and checks this equality for $n=1$ only.
But the operators $T,L$ under the consideration do not commute, and easy examples show that  Gorbatyuk's equality  (\ref{eqGo})
is not valid even for $n =2$. Nevertheless, it turns out that the asymptotic equation (\ref{eqGoT}) is true.
In the present paper we prove a stronger statement.

\begin{thm}\label{thm2}
When $n\to \infty$ we have
$$
 \widehat{m}(2^n)=1+O\left(\frac{1}{\log n}\right).
$$
\end{thm}

Our proof of this Theorem \ref{thm2} is based on the result by  Kesseböhmer and  Stratmann (\ref{kss}).

\subsection{Structure of the paper}
The rest of the paper is organised as follows.

In Section  \ref{GK} we give a proof of Theorem \ref{thmBV} and deduce Corollaries 1, 2, 3.
In particular, to do this in Subsection \ref{sss2} we deal with  the simplest properties of operator $T$ which we use in all of our proofs.
In Section \ref{simp} we 
express the reminder $R_n$ by means of the auxiliary function 
$\rho(x)=\{x\}-\frac{1}{2}$ and
deduce from Theorem \ref{thmBV} equality (\ref{simpler}).
Section \ref{uuu}
is devoted to more detailed analysis of  action powers  of operator $T$ and to the proof of Theorem \ref{thm1}.
In Section \ref{vvv} we prove Theorem \ref{thm2}.

\section{ General construction}\label{GK}

In Subsection \ref{sss1}  below we reduce the problem of obtaining an upper bound for the reminder $R_n$ to the problem of estimating of a certain integral
$\beta_n$. In Subsection \ref{sss2} we introduce operator $T$ and study its properties. In Subsection \ref{sss3} we finalise the proof of (\ref{simpler}).

\subsection{Operator $T$: the simplest properties}\label{sss2}

 We consider a linear operator $T$ on the space of continuous functions
 $C\left[\frac{1}{2},1\right]$ defined by formula
 \begin{equation}\label{ope}
 T\, 
 f(x) = \frac{f\left(\frac{1}{2-x}\right)}{(2-x)^2}-\frac{f\left(\frac{1}{1+x}\right)}{(1+x)^2}.
 \end{equation}

First of all we observe that operator $T$ defined in (\ref{ope}) has the following obvious properties.
 
 \vskip+0.3cm
 
 {\bf Property 1.} \,\, $T$ maps (strictly) increasing functions to (strictly) increasing functions.
 
 \vskip+0.3cm
{\bf Property 2.} \,\, $\forall f\in C\left[\frac{1}{2},1\right]$, $(Tf)\left(\frac{1}{2}\right)=0$. And as a corollary, by Property 1, $Tf$ is nonnegative whenever $f$ is increasing.\

\vskip+0.3cm
{\bf Property 3.} \,\, $\forall k\ge1$, $(T^kf)(1)=(Tf)(1)$. And as a corollary, by Property 1, $T^kf$ is bounded uniformly in $k$ if $f\in C\left[\frac{1}{2},1\right]$ is increasing.

\vskip+0.3cm

 Then we prove two lemmas.
 
\begin{lem}\label{l3}
Let $j(x)$  be 
 $1$-periodic
 odd piecewise continuous bounded function,
  and  $f_0(x)$ be continuous function  defined on $\left[\frac{1}{2},1\right]$.
 Let 
  $f_n=T^nf_0$.
  Then
  \begin{equation}\label{media}
  \int_{\frac{1}{2}}^1j(2^n?(x))f_0(x)\text{d}x=\int_{\frac{1}{2}}^1j(?(x))f_n(x)\text{d}x.
  \end{equation}
 \end{lem}
\noindent
{\bf Proof.} 
We should note that for any  continuous function $g(x)$ the composition $j(g(x))$ will be a Riemann integrable function.
For all $0\le k\le n-1$
we have
$$
\begin{aligned}&\int_{\frac{1}{2}}^1j(2^{n-k}?(x))f_k(x)\text{d}x\xlongequal{x:=\frac{1}{1+t}}\int_1^0j\left(2^{n-k}?\left(\frac{1}{1+t}\right)\right)f_k\left(\frac{1}{1+t}\right)(-(1+t)^{-2})\text{d}t\\=&\int_0^1j\left(2^{n-k}\left(1-?\left(\frac{t}{1+t}\right)\right)\right) f_k\left(\frac{1}{1+t}\right)(1+t)^{-2}\text{d}t=\int_0^1j\left(2^{n-k}\left(1-\frac{?(t)}{2}\right)\right)f_k\left(\frac{1}{1+t}\right)(1+t)^{-2}\text{d}t\\=&-\int_0^1j(2^{n-k-1}?(t))f_k\left(\frac{1}{1+t}\right)(1+t)^{-2}\text{d}t=-\left(\int_0^{\frac{1}{2}}+\int_{\frac{1}{2}}^1\right).
\end{aligned}
$$
Then we transform the first integral in the right hand side as
$$
\begin{aligned}&\int_0^{\frac{1}{2}}j(2^{n-k-1}?(t))f_{k}\left(\frac{1}{1+t}\right)(1+t)^{-2}\text{d}t\xlongequal{x:=1-t}-\int_1^{\frac{1}{2}}j(2^{n-k-1}?(1-x))f_k\left(\frac{1}{2-x}\right)(2-x)^{-2}\text{d}x\\=&\int_{\frac{1}{2}}^1j(2^{n-k-1}(1-?(x)))f_k\left(\frac{1}{2-x}\right)(2-x)^{-2}\text{d}x=-\int_{\frac{1}{2}}^1j(2^{n-k-1}?(x))f_k\left(\frac{1}{2-x}\right)(2-x)^{-2}\text{d}x.\end{aligned}
$$
Now we continue with the sum of integrals $\int_0^{\frac{1}{2}}+\int_{\frac{1}{2}}^1$
as follows:
$$\int_{\frac{1}{2}}^1j(2^{n-k}?(x))f_k(x)\text{d}x=\int_{\frac{1}{2}}^1j(2^{n-k-1}?(x))\left(\frac{f_k\left(\frac{1}{2-x}\right)}{(2-x)^2}-\frac{f_k\left(\frac{1}{1+x}\right)}{(1+x)^2}\right)\text{d}x=\int_{\frac{1}{2}}^1j(2^{n-k-1}?(x))f_{k+1}(x)\text{d}x.$$
Lemma is proven.$\Box$
\vskip+0.3cm

 The next lemma provides a recursive  equality  for the values of functions $f_n = T^nf_0$ for function $f_0(x)$ defined on $\left[\frac{1}{2},1\right]$.
 
 \begin{lem}\label{l50}
  Let  $ f_0 (x) = f(x)$ and
  $ f_n (x)  = (T^n f_0)(x)$. Let  $ y\ge 1$.
  Then 
 \begin{equation}\label{recurrence}
 f_n
 \left(
  \frac{y}{y+1}\right) =
  (y+1)^2 \left(
  \frac{f_0
  \left(
  \frac{y+n}{y+n+1}
  \right)}{
  (y+n+1)^2}
   -\sum_{k=1}^n 
   \frac{f_{n-k} 
   \left(
  \frac{y+k}{2y+2k-1}
  \right)}
  {(2y+2k-1)^2}
  \right).
  \end{equation}
  \end{lem}

    \noindent
{\bf Proof.} 
We  proceed by induction. By the definition (\ref{ope}) we have
$$
f_1 \left(
  \frac{y}{y+1}\right) 
  =
  \left(
  \frac{y+1}{y+2}\right)^2 
  f_0  \left(
  \frac{y+1}{y+2}\right) -
    \left(
  \frac{y+1}{2y+1}\right)^2 
  f_0  \left(
  \frac{y+1}{2y+1}\right) 
 ,
  $$
 and this proves the base of induction for  $n=1$.
 Now we proceed with the step. Assume (\ref{recurrence}) we get
  $$
  f_n
 \left(
  \frac{y+1}{y+2}\right) =
  (y+2)^2 \left(
  \frac{f_0
  \left(
  \frac{y+n+1}{y+n+2}
  \right)}{
  (y+n+2)^2}
   -\sum_{k=1}^n 
   \frac{f_{n-k} 
   \left(
  \frac{y+k+1}{2y+2k+1}
  \right)}
  {(2y+2k+1)^2}
  \right).
  $$
 Now applying
 $f_{n+1} (x) = (Tf_n)(x)$ we continue with
 $$
  f_{n+1}
 \left(
  \frac{y}{y+1}\right) 
    =
  \left(
  \frac{y+1}{y+2}\right)^2 
  f_n  \left(
  \frac{y+1}{y+2}\right) -
    \left(
  \frac{y+1}{2y+1}\right)^2 
  f_n  \left(
  \frac{y+1}{2y+1}\right) 
 $$
 $$
    =  (y+1)^2 \left(
  \frac{f_0
  \left(
  \frac{y+n+1}{y+n+2}
  \right)}{
  (y+n+2)^2}
   -
   \sum_{k=1}^n 
   \frac{f_{n-k} 
   \left(
  \frac{y+k+1}{2y+2k+1}
  \right)}
  {(2y+2k+1)^2}
  \right)
  -
    \left(
  \frac{y+1}{2y+1}\right)^2 
  f_n  \left(
  \frac{y+1}{2y+1}\right) 
 $$
 $$
 =
 (y+1)^2 \left(
  \frac{f_0
  \left(
  \frac{y+n+1}{y+n+2}
  \right)}{
  (y+n+2)^2}
   -
   \sum_{k=1}^{n+1} 
   \frac{f_{n+1-k} 
   \left(
  \frac{y+k}{2y+2k-1}
  \right)}
  {(2y+2k-1)^2}
  \right).
 $$
 Everything is proven.$\Box$

\vskip+0.3 cm

The following statement is an immediate corollary of Lemma \ref{l50}.

 \begin{lem}\label{l5} Supppose that$$f_0(x)\in C\left[\frac{1}{2},1\right]$$has bounded variation. Therefore we have $f_0(x)=g_0(x)-h_0(x)$, where $g_0(x)$ and $h_0(x)$ are nonnegative bounded increasing functions. Let$$M:=\max\left(\max_{x\in\left[\frac{1}{2},1\right]}g_0(x),\max_{x\in\left[\frac{1}{2},1\right]}h_0(x)\right).$$
  Then for any $ n\ge 0$  and $ y\ge 1$ we have inequality
  \begin{equation}\label{l5_1}
   \left|f_n\left(\frac{y}{y+1}\right)\right|\le2M\left(\frac{y+1}{y+n+1}\right)^2, 
   \end{equation}
   and moreover, for any $ n\ge k \ge 1$ and $y\ge 1$ one has
    \begin{equation}\label{l5_2}
   \left|f_{n-k}\left(\frac{y+k}{2y+2k-1}\right)\right|\le\frac{32M}{(n-k)^2}. 
   \end{equation}

  \end{lem}

    \noindent
{\bf Proof.}  
{Inequality}
(\ref{l5_1}) with $f_0$ replaced by $g_0$ and $h_0$ follows immediately from {\bf Properties 2}, {\bf 3} and Lemma \ref{l50} by ignoring the negative summand in the right hand side of (\ref{recurrence}).

To obtain (\ref{l5_2}) we use (\ref{l5_1}) with $n_1= n-k$ instead of $n$ and  $y_1= \frac{y+k}{y+k-1}\ge 1$ instead of $y$:
$$
   \left|f_{n-k}
   \left(\frac{y+k}{2y+2k-1}\right)\right|
   =
      \left|f_{n_1}
   \left(\frac{y_1}{y_1+1}\right)\right|
   \le
  2M\left(\frac{y_1+1}{y_1+n_1+1}\right)^2=
    2M\left(\frac{2y+2k-1}{(n-k+2)(y+k-1)+1}\right)^2
    \le \frac{32M}{(n-k)^2},
$$
as $ y, k \ge 1$.
Everything is proven.
$\Box$

\subsection{Proof of Theorem \ref{thmBV}} \label{aaa}

In order to prove Theorem  \ref{thmBV} 
we need to use only inequality (\ref{l5_1}) of Lemma \ref{l5}.

 For  $x\in\left[\frac{1}{2},1\right)$ we use the identity
$$
x= \frac{y}{y+1},\,\,\,\,\text{where}\,\,\,\,  y = y(x)  = \frac{x}{1-x}.
$$
{
We use Lemma \ref{l3} and
 estimate the right hand side from (\ref{media}) by means of  (\ref{l5_1}). So we get}
 $$
  \int_{\frac{1}{2}}^1j(2^n?(x))F_0(x)\text{d}x\le
M_1\int_{\frac{1}{2}}^1\left( \frac{y+1}{y+n+1} \right)^2{\rm d}x=
   M_1 \int_1^\infty\frac{\text{d}y}{(y+n+1)^2}=\frac{M_1}{n+2} = O\left(
  \frac{1}{n}\right),
 $$
 with
 $
 M_1 =
   2M \max_{x\in \mathbb{R}}  |j(x)|.
   $
We proved (\ref{oo}).$\Box$

\subsection{Proof of Corollaries 1, 2, 3}\label{bbb}

All the corollaries easily follow from Theorem \ref{thmBV}.

\vskip+0.3cm
\noindent
{\bf Proof of Corollary 1}.
 In Theorem \ref{thmBV} we put 
 $$
 j(x)=\begin{cases}
 (-1)^{[2x]},&x\notin\frac{1}{2}\mathbb{Z};\\0,&x\in\frac{1}{2}\mathbb{Z}.
 \end{cases}
 $$ 
 As segments $[\xi_{k,n+1},\xi_{k+1,n+1}]$ form a partition of the interval $\left[\frac{1}{2},1\right]$,
 we have
 $$\int_{\frac{1}{2}}^1j(2^n?(x))F(x)\text{d}x=\sum_{k=2^n}^{2^{n+1}-1}\int_{\xi_{k,n+1}}^{\xi_{k+1,n+1}}j(2^n?(x))F(x)dx=\sum_{k=2^n}^{2^{n+1}-1}(-1)^k\int_{\xi_{k,n+1}}^{\xi_{k+1,n+1}}F(x)dx,$$
 and by Theorem \ref{thmBV} we are done. $\Box$

\vskip+0.3cm
\noindent
{\bf Proof of Corollary 2}.
For $m$ positive by taking $F(x)=x^{m-1}$ 
in Corollary 1, 
we obtain
$$
\xi_{2^{n-1},n}^m+\xi_{2^n,n}^m
-2\sum_{k=2^{n-1}}^{2^n}(-1)^k\xi_{k,n}^m=
\sum_{k=2^{n-1}}^{2^n-1}(-1)^k(\xi_{k+1,n}^m-\xi_{k,n}^m)=  O(n^{-1}).
$$
But  $\xi_{2^{n-1},n}=\frac{1}{2}$ and $\xi_{2^n,n}=1$,
so 
$$
 \sum_{k=2^{n-1}}^{2^n}(-1)^k\xi_{k,n}^m=\frac{1}{2}+\frac{1}{2^{m+1}}
 + O(n^{-1}), 
$$
and
Corollary 2 follows.$\Box$

\vskip+0.3cm
\noindent
{\bf Proof of Corollary 3}.
  For positive $m$
  we take into account equalities 
  { (\ref{sYM})}
  which together with Corollary 2
  give us
 $$
 \begin{aligned}\sum_{k=0}^{2^n}(-1)^k\xi_{k,n}^m&=\sum_{k=2^{n-1}}^{2^n}(-1)^k\xi_{k,n}^m+\sum_{k=0}^{2^{n-1}}(-1)^k\xi_{k,n}^m-\xi_{2^{n-1},n}^m\\&=\sum_{k=2^{n-1}}^{2^n}(-1)^k\xi_{k,n}^m+\sum_{k=2^{n-1}}^{2^n}(-1)^k(1-\xi_{k,n})^m-\frac{1}{2^m}\\&=\sum_{k=2^{n-1}}^{2^n}(-1)^k\xi_{k,n}^m+\sum_{l=0}^m\binom{m}{l}(-1)^l\sum_{k=2^{n-1}}^{2^n}(-1)^k\xi_{k,n}^l-\frac{1}{2^m}\\&=\frac{1}{2}+\frac{1}{2^{m+1}}+\sum_{l=0}^m\binom{m}{l}(-1)^l\left(\frac{1}{2}+\frac{1}{2^{l+1}}\right)-\frac{1}{2^m}
+ O(n^{-1})
 =\frac{1}{2} +
 O(n^{-1}),\,\,n\to\infty,
 \end{aligned},
 $$
 and everything is proven.$\Box$

\section{Proof of  simpler equality (\ref{simpler})}\label{simp}

As we have mentioned before, equality  (\ref{simpler}) can be deduced directly from Theorem \ref{thmBV}.

\subsection{ Function $\rho(x)=\{x\}-\frac{1}{2}$: auxiliary results }\label{sss1}

Here we prove several auxiliary statements.
Let $\rho(x)=\{x\}-\frac{1}{2}$, where $\{x\}=x-\lfloor x\rfloor$ is the fractional part of $x$ 

\begin{lem}\label{l1}
For the remainder $R_n$ we have identity
\begin{equation}\label{rem}
R_n=
 \sum_{j=1}^{2^n}\left(\xi_{j,n}-\frac{j}{2^n}\right)^2-
 2^n\int_0^1(?(x)-x)^2\text{d}x=\int_0^1\rho(2^nx)\text{\rm d}(m(x)-x)^2.
 \end{equation}
\end{lem}

\noindent
{\bf Proof.}
First we observe that integration by parts gives 
$$\int_0^1(?(x)-x)^2\text{d}x=\int_0^1(m(x)-x)^2\text{d}x.$$
Then
$$
\begin{aligned}
&2^n\int_0^1(?(x)-x)^2\text{d}x-\sum_{j=1}^{2^n}\left(\xi_{j,n}-\frac{j}{2^n}\right)^2=2^n\int_0^1(m(x)-x)^2\text{d}x-\sum_{j=1}^{2^n}\left(\xi_{j,n}-\frac{j}{2^n}\right)^2\\=&2^n\int_0^1(m(x)-x)^2\text{d}x-\sum_{j=1}^{2^n}\left(m\left(\frac{j}{2^n}\right)-\frac{j}{2^n}\right)^2=2^n\sum_{j=1}^{2^n}\left(\int_{\frac{j-1}{2^n}}^{\frac{j}{2^n}}(m(x)-x)^2\text{d}x-\frac{1}{2^n}\left(m\left(\frac{j}{2^n}\right)-\frac{j}{2^n}\right)^2\right).
\end{aligned}
$$
Integrating by part we obtain
$$
\begin{aligned}
&2^n\sum_{j=1}^{2^n}\left(\int_{\frac{j-1}{2^n}}^{\frac{j}{2^n}}(m(x)-x)^2\text{d}x-\frac{1}{2^n}\left(m\left(\frac{j}{2^n}\right)-\frac{j}{2^n}\right)^2\right)=-2^n\sum_{j=1}^{2^n}\int_{\frac{j-1}{2^n}}^{\frac{j}{2^n}}\left(x-\frac{j-1}{2^n}\right)\text{d}(m(x)-x)^2\\=&-2^n\sum_{j=1}^{2^n}\int_{\frac{j-1}{2^n}}^{\frac{j}{2^n}}\left(x-\frac{j-1}{2^n}-\frac{1}{2^{n+1}}\right)\text{d}(m(x)-x)^2=-\int_0^1\rho(2^nx)\text{d}(m(x)-x)^2.
\end{aligned}
$$
Lemma is proven.$\Box$

\vskip+0.3cm
Now the integral in the right hand side of (\ref{rem}) we represent as a sum
\begin{equation}\label{sum1}
\int_0^1\rho(2^nx)\text{d}(m(x)-x)^2=2\int_0^1\rho(2^nx)(m(x)-x)\text{d}m(x)+2\int_0^1\rho(2^nx)(x-m(x))\text{d}x.
\end{equation}
The next lemma deals with the second summand from (\ref{sum1}).

\begin{lem} \label{l2}
For the second summand in  (\ref{sum1}) we have
 $$\int_0^1\rho(2^nx)(x-m(x))\text{\rm d}x=O\left(\frac{1}{2^n}\right).$$
\end{lem}

\noindent
{\bf Proof.}
We rewrite the integral as
$$
\int_0^1\rho(2^nx)(x-m(x))\text{d}x=\sum_{j=1}^{2^n}\int_{\frac{j-1}{2^n}}^{\frac{j}{2^n}}\rho(2^nx)(x-m(x))\text{d}x\xlongequal{u:=2^nx}\frac{1}{2^n}\sum_{j=1}^{2^n}\int_{j-1}^j\rho(u)\left(\frac{u}{2^n}-m\left(\frac{u}{2^n}\right)\right)\text{d}u.
$$
To continue with the last integral we observe that 
$$
\sum_{j=1}^{2^n}\int_{j-1}^j\rho(u)\left(\frac{j}{2^n}-m\left(\frac{j}{2^n}\right)\right)\text{d}u=\sum_{j=1}^{2^n}\left(\int_{j-1}^j\rho(u)\text{d}u\right)\cdot\left(\frac{j}{2^n}-m\left(\frac{j}{2^n}\right)\right)=0.
$$
Now
$$
\begin{aligned}
&\left|\frac{1}{2^n}\sum_{j=1}^{2^n}\int_{j-1}^j\rho(u)\left(\frac{u}{2^n}-m\left(\frac{u}{2^n}\right)\right)\text{d}u\right|=\left|\frac{1}{2^n}\sum_{j=1}^{2^n}\int_{j-1}^j\rho(u)\left(\frac{u}{2^n}-\frac{j}{2^n}-m\left(\frac{u}{2^n}\right)+m\left(\frac{j}{2^n}\right)\right)\text{d}u\right|\\\le&\frac{1}{2^n}\sum_{j=1}^{2^n}\int_{j-1}^j|\rho(u)|\left(\left|\frac{u}{2^n}-\frac{j}{2^n}\right|+\left|-m\left(\frac{u}{2^n}\right)+m\left(\frac{j}{2^n}\right)\right|\right)\text{d}u\le\frac{1}{2^n}\sum_{j=1}^{2^n}\frac{1}{2}\left(\frac{j}{2^n}-\frac{j-1}{2^n}+\xi_{j,n}-\xi_{j-1,n}\right)=\frac{1}{2^n},\end{aligned}
$$
and
lemma is proven.$\Box$

\subsection{End of the proof of  (\ref{simpler})} \label{sss3}
From Lemmas \ref{l1}, \ref{l2} and (\ref{sum1}) it follows that for the remainder $R_n$ we have equality
\begin{equation}\label{R_beta}
R_n = 2\beta_n + O(2^{-n}),
\end{equation}
where
$$
 \beta_n = \int_0^1\rho(2^nx)(m(x)-x)\text{d}m(x)=\int_0^1\rho(2^n?(x))(x-?(x))\text{d}x.
$$
We denote $ f_0(x) = x-?(x)$  and notice that $f_0(x)=-f_0(1-x)$.
As $\rho(kx),\,k\in\mathbb{Z}$ is an odd function of period 1, we have
$$
\beta_n=
\int_0^{\frac{1}{2}}\rho(2^n?(x))f_0(x)\text{d}x\xlongequal{u=1-x}\int_{\frac{1}{2}}^1\rho(2^n?(u))f_0(u)\text{d}u.
$$
Now we apply Theorem \ref{thmBV}
for functions $j(x) = \rho(x)$
and 
$F(x) = f_0(x) $
to get
$$
\beta_n  = O(n^{-1}),\,\,\,n\to \infty.
$$
By (\ref{R_beta}), this gives (\ref{simpler}).$\Box$

 \section{Proof of Theorem \ref{thm1}}\label{uuu}
  In this section let $f_0(x):=x-?(x)$.\\
  The idea of the proof is that by Lemma \ref{l50} and the fact  that $f_0\left(\frac{y+n+1}{y+n+2}\right)=O\left(\frac{1}{y+n}\right)$, what is left 
  {is to estimate the sum}
 $$ \sum_{k=1}^{n} 
   \frac{f_{n-k} 
   \left(
  \frac{y+k}{2y+2k-1}
  \right)}
  {(2y+2k-1)^2}.$$
  We estimate this sum by (\ref{l5_2}) when $k$ is not very close to $n$, while when $k$ is close to $n$ we need the following statement.

  \begin{lem}\label{l6v} { {\bf (Main Lemma)}}
  {
  For all   $ n\ge 0$ we have
$$f_n(x)\ll\frac{1}{\log n}\left( x-\frac{1}{2}\right),\,\,\,
\text{
uniformly on}\,\,\,x 
\in\left[\frac{1}{2},\frac{2}{3}\right].
$$}
\end{lem}
 
{
\subsection{Proof of Main Lemma}
}

  To prove this lemma we write $f_0(x)=x-\frac{1}{2}-\left(?(x)-\frac{1}{2}\right)$ and deal with
  { functions}
   $x$, $\frac{1}{2}$, $?(x)-\frac{1}{2}$ separately in Lemmas \ref{l6}, \ref{hamburg} and \ref{l8} 
   {below and then use }
    the linearity of $T$.

\begin{lem}\label{l6}  
{
We have
$$T^nx\ll\frac{1}{\log n}\left( x-\frac{1}{2}\right),\,\,\,
\text{uniformly in}\,\,\,x \in\left[\frac{1}{2},\frac{2}{3}\right].$$
}
\end{lem}

\noindent
{\bf Proof.}
  { First of all we prove by induction that 
  that for every positive integer $n\ge2$ one has}
\begin{equation}\label{bayern}
T^{n}x=\sum_{\frac{p}{q}\in W_{n-1}}\left((q-px)\left(\frac{1}{(q_+-p_+x)^3}-\frac{1}{(q_--p_-x)^3}\right)-(q-p(1-x))\left(\frac{1}{(q_+-p_+(1-x))^3}-\frac{1}{(q_--p_-(1-x))^3}\right)\right),
\end{equation}
where $\frac{p_-}{q_-}<\frac{p}{q}<\frac{p_+}{q_+}$ are the two descendants of $\frac{p}{q} \in W_{n-1}$ in $W_n$.\\
{ As for the base of induction, one can check formula  (\ref{bayern}) for $n=2$ by straightforward calculation}.\\
  Observe that since $\frac{1}{2-x}$ is strictly increasing and $\frac{1-x}{2-x}$ is strictly decreasing, 
 {we have}
  $\frac{q_-}{2q_--p_-}<\frac{q}{2q-p}<\frac{q_+}{2q_+-p_+}$
  and $\frac{q_+-p_+}{2q_+-p_+}<\frac{q-p}{2q-p}<\frac{q_--p_-}{2q_--p_-}$.
 { Recall that we have the rule (\ref{wn+1}) how to get elements of $W_{n+1}$ from the elements of $W_n$.
 So 
  $\frac{q_-}{2q_--p_-},\frac{q_+}{2q_+-p_+} \in W_{n+1}$
   are exactly the descendants of $\frac{q}{2q-p}\in W_n$.
   At the same time fractions
   $\frac{q_+-p_+}{2q_+-p_+},\frac{q_--p_-}{2q_--p_-}\in W_{n+1}$
      are exactly the descendants of $
   \frac{q-p}{2q-p}\in W_n$.
  } 
  
  Then for $n+1$
  {we have}
  $$\begin{aligned}T^{n+1}x&=T(T^nx)\\&=\sum_{\frac{p}{q}\in W_{n-1}}{\left(((2q-p)-qx)\left(\frac{1}{((2q_+-p_+)-q_+x)^3}-\frac{1}{((2q_--p_-)-q_-x)^3}\right)\right.}\\&{\left.\quad\quad\quad-((2q-p)-(q-p)x)\left(\frac{1}{((2q_+-p_+)-(q_+-p_+)x)^3}-\frac{1}{((2q_--p_-)-(q_--p_-)x)^3}\right)\right)}\\&\quad-\sum_{\frac{p}{q}\in W_{n-1}}{\left(((2q-p)-q(1-x))\left(\frac{1}{((2q_+-p_+)-q_+(1-x))^3}-\frac{1}{((2q_--p_-)-q_-(1-x))^3}\right)\right.}\\&{\left.\quad\quad\quad-((2q-p)-(q-p)(1-x))\left(\frac{1}{((2q_+-p_+)-(q_+-p_+)(1-x))^3}-\frac{1}{((2q_--p_-)-(q_--p_-)(1-x))^3}\right)\right)}\\&=\sum_{\frac{p}{q}\in W_{n-1}}{\left(((2q-p)-qx)\left(\frac{1}{((2q_+-p_+)-q_+x)^3}-\frac{1}{((2q_--p_-)-q_-x)^3}\right)\right.}\\&{\left.\quad\quad\quad-((2q-p)-q(1-x))\left(\frac{1}{((2q_+-p_+)-q_+(1-x))^3}-\frac{1}{((2q_--p_-)-q_-(1-x))^3}\right)\right)}\\&\quad+\sum_{\frac{p}{q}\in W_{n-1}}{\left(((2q-p)-(q-p)x)\left(\frac{1}{((2q_--p_-)-(q_--p_-)x)^3}-\frac{1}{((2q_+-p_+)-(q_+-p_+)x)^3}\right)\right.}\\&{\left.\quad\quad\quad-((2q-p)-(q-p)(1-x))\left(\frac{1}{((2q_--p_-)-(q_--p_-)(1-x))^3}-\frac{1}{((2q_+-p_+)-(q_+-p_+)(1-x))^3}\right)\right)}\\&=\sum_{\frac{p}{q}\in W_n}\left((q-px)\left(\frac{1}{(q_+-p_+x)^3}-\frac{1}{(q_--p_-x)^3}\right)-(q-p(1-x))\left(\frac{1}{(q_+-p_+(1-x))^3}-\frac{1}{(q_--p_-(1-x))^3}\right)\right).
\end{aligned}$$
  Hence by induction (\ref{bayern}) is proven.
  
  {
 Now for the function $u(x) = T^n(x)$ we use mean value theorem and get $\xi \in \left[\frac{1}{2},\frac{2}{3}\right]$ such that 
 $$
 u(x) - u\left(\frac{1}{2}\right) = \left(x - \frac{1}{2}\right) u'(\xi).
 $$
 Recall that by {\bf  Property 2} $ u\left(\frac{1}{2}\right)  =0$ So from  (\ref{bayern}) for $x\in\left[\frac{1}{2},\frac{2}{3}\right]$ we get
   $$
   T^nx=\left(x - \frac{1}{2}\right)
  (
   \sum_{\frac{p}{q}\in W_{n-1}}\left(-p\left(\frac{1}{(q_+-p_+\xi)^3}-\frac{1}{(q_--p_-\xi)^3}\right)+(q-p\xi)\left(\frac{3p_+}{(q_+-p_+\xi)^4}-\frac{3p_-}{(q_--p_-\xi)^4}\right)\right)
   $$
  $$
 + \sum_{\frac{p}{q}\in W_{n-1}}\left(-p\left(\frac{1}{(q_+-p_+\xi_1)^3}-\frac{1}{(q_- -p_- \xi_1)^3}\right)+(q-p\xi_1)\left(\frac{3p_+}{(q_+-p_+\xi_1)^4}-\frac{3p_-}{(q_--p_-\xi_1)^4}\right)\right) ), .
   $$
  {where} $ \xi_1 = 1- \xi $.}
   Hence
   $$
   |T^nx|\le2\left(x-\frac{1}{2}\right)\sum_{\frac{p}{q}\in W_{n-1}}\left({\frac{27p}{q_+^3}}+\frac{27p}{q_-^3}+\frac{243qp_+}{q_+^4}+\frac{243qp_-}{q_-^4}\right)\le
   {1080}
   \left(x-\frac{1}{2}\right)\sum_{\frac{p}{q}\in W_n}\frac{1}{q^2}\ll\frac{1}{\log n}\left(x-\frac{1}{2}\right).
   $$
Here we used (\ref{kss}) from \cite{ks} to deduce the coefficient $\frac{1}{\log n}$.$\Box$\\
\begin{lem}\label{hamburg}
  {We have$$T^n1\ll\frac{1}{\log n}\left(x-\frac{1}{2}\right),\,\,\,\text{uniformly in}\,\,\, x\in\left[\frac{1}{2},\frac{2}{3}\right]$$
  }
\end{lem}
\noindent
{\bf Proof.}
  Analogous to Lemma \ref{l6}, we first prove by induction that for all $n\ge1$, holds
\begin{equation}\label{frankfurt}
T^{n}1=\sum_{\frac{p}{q}\in W_n}\varepsilon_{\frac{p}{q}}\left(\frac{1}{(q-px)^2}-\frac{1}{(q-p(1-x))^2}\right),
\end{equation}
where $\varepsilon_\frac{p}{q}\in\{-1,1\}$.\\
  The $n=1$ case is easy to check. Suppose for $n$ (\ref{frankfurt}) holds. Then for $n+1$,$$\begin{aligned}T^{n+1}1&=T(T^n1)\\&=\sum_{\frac{p}{q}\in W_n}\varepsilon_{\frac{p}{q}}\left(\frac{1}{((2q-p)-qx)^2}-\frac{1}{((2q-p)-(q-p)x)^2}\right)-\sum_{\frac{p}{q}\in W_n}\varepsilon_{\frac{p}{q}}\left(\frac{1}{((q-p)+qx)^2}-\frac{1}{(q+(q-p)x)^2}\right)\\&=\sum_{\frac{p}{q}\in W_n}\varepsilon_{\frac{p}{q}}\left(\frac{1}{((2q-p)-qx)^2}-\frac{1}{((2q-p)-q(1-x))^2}\right)\\&\quad+\sum_{\frac{p}{q}\in W_n}(-\varepsilon_{\frac{p}{q}})\left(\frac{1}{((2q-p)-(q-p)x)^2}-\frac{1}{((2q-p)-(q-p)(1-x))^2}\right)\\&=\sum_{\frac{p}{q}\in W_{n+1}}\varepsilon_{\frac{p}{q}}\left(\frac{1}{(q-px)^2}-\frac{1}{(q-p(1-x))^2}\right).\end{aligned}$$
  Hence by induction (\ref{frankfurt}) is proven.
  
  {
  Now  as in the end of the proof of  Lemma  \ref{l6} we apply mean-value theorem.}$\Box$
 \vskip+0.3cm 
  
  In order to continue with $?(x)-\frac{1}{2}$ we need one more auxiliary statement.
We should consider one more operator $T_+ :C\left[\frac{1}{2},1\right]\to C\left[\frac{1}{2},1\right]$ defined by
$$
T_+f(x)=\frac{f\left(\frac{1}{2-x}\right)}{(2-x)^2}+\frac{f\left(\frac{1}{1+x}\right)}{(1+x)^2}.
$$
This operator was used by Golubeva in her paper \cite{g}.
\vskip+0.3cm

We need to consider one more family of functions.
Let 
$$
\psi_0(x)=1\,\,\,\,\,\text{ and} \,\,\,\,\, \psi_n(x)=(T_+)^n\psi_0 (x),
$$
Applying (\ref{wn+1}) it is not difficult to obtain
\begin{equation}
\label{eqn-30}\psi_n(x)=\sum_{\frac{p}{q}\in W_n}\left(\frac{1}{(q-px)^2}+\frac{1}{(q-p(1-x))^2}\right)
\end{equation}
 {
(formula (\ref{eqn-30}) for $\psi_n(x)$, in particular, is given in  Lemma 2} from \cite{g}).

\begin{lem}\label{l8} 
  {We have
  $$T^n\left(?(x)-\frac{1}{2}\right)\ll\frac{1}{\log n}\left(x-\frac{1}{2}\right),\,\,\, \text{uniformly in}\,\,\,
 x\in\left[\frac{1}{2},\frac{2}{3}\right].
  $$}
\end{lem}

\noindent
{\bf Proof.} 
  We first prove that $?(x)-\frac{1}{2}\le201\left(x-\frac{1}{2}\right)$. This inequality obviously holds for $x\in\left(\frac{101}{201},\frac{2}{3}\right]$, since in this case $201\left(x-\frac{1}{2}\right)>1$.\\For $x\in\left[\frac{1}{2},\frac{101}{201}\right]$, $x=[0;1,1,a_3,\ldots]$ with $a_3\ge 5$. 
Hence we have $x-\frac{1}{2}\ge\frac{1}{4a_3+8}$ (this can be seen from, e.g.  Perron's formula). 
By the explicit formula (\ref{eqn-1}) for $?(x)$ we have$$?(x)-\frac{1}{2}\le\frac{1}{2^{a_3+1}}<\frac{1}{4a_3+8} \le  x-\frac{1}{2}.$$Thus $?(x)-\frac{1}{2}\le201\left(x-\frac{1}{2}\right)$, on $\left[\frac{1}{2},\frac{2}{3}\right]$.\\
  Next we prove an identity analogous to (\ref{bayern}) and (\ref{frankfurt}). We have, for all $x\in\left[\frac{1}{2},\frac{2}{3}\right]$, $\forall n\ge1$,
\begin{equation}\label{dortmund}
T^n\left(?(x)-\frac{1}{2}\right)=\frac{1}{2^n}\left(?(x)-\frac{1}{2}\right)\psi_n(x)+\sum_{\frac{p}{q}\in W_n}{c_{\frac{p}{q}}}\left(\frac{1}{(q-px)^2}-\frac{1}{(q-p(1-x))^2}\right).
\end{equation}

The values of $ c_{\frac{p}{q}}$ for $\frac{p}{q}\in W_n$ satisfy inequality $|c_{\frac{p}{q}}|\le\frac{1}{2}-\frac{1}{2^{n+1}}.$
  Taking into account equalities
$$?\left(\frac{1}{2-x}\right) - \frac{1}{2} = \frac{?(x)}{2},\,\,\,\,\,
?\left(\frac{1}{1+x}\right) - \frac{1}{2} = \frac{1}{2} - \frac{?(x)}{2},
$$
which follow from (\ref{eqn-2}), we proceed by induction. For $n=1$ the set $W_1$ contains just one number $\frac{1}{2}$, namely $p=1,q=2$, and (\ref{dortmund}) follows from (\ref{psi1}) 
with $c_{\frac{1}{2}} =\frac{1}{4}$ as
$$
h_1(x) =
\frac{1}{2} \left(
?(x) - \frac{1}{2}\right) 
\left(
\frac{1}{(2-x)^2}+\frac{1}{(1+x)^2}\right) +
\frac{1}{4} \left(
\frac{1}{(2-x)^2}-\frac{1}{(1+x)^2}\right).
$$

  Suppose (\ref{dortmund}) holds for $n$, then for $n+1$,$$
\begin{aligned}h_{n+1}(x)&=T\left(\frac{1}{2^n}\left(?(x)-\frac{1}{2}\right)\psi_n(x)+\sum_{\frac{p}{q}\in W_n}c_{\frac{p}{q}}\left(\frac{1}{(q-px)^2}-\frac{1}{(q-p(1-x))^2}\right)\right)\\&=\frac{1}{2^n}\left(\frac{\left(?\left(\frac{1}{2-x}\right)-\frac{1}{2}\right)\psi_n\left(\frac{1}{2-x}\right)}{(2-x)^2}-\frac{\left(?\left(\frac{1}{1+x}\right)-\frac{1}{2}\right)\psi_n\left(\frac{1}{1+x}\right)}{(1+x)^2}\right)\\&\quad+\sum_{\frac{p}{q}\in W_n}\left(c_{\frac{p}{q}}\left(\frac{1}{((2q-p)-qx)^2}-\frac{1}{((q-p)+qx)^2}\right)-c_{\frac{p}{q}}\left(\frac{1}{((2q-p)-(q-p)x)^2}-\frac{1}{(q+(q-p)x)^2}\right)\right).
\end{aligned}
$$ 
We continue with the first summand from the right hand side here.
$$
\frac{1}{2^n}\left(\frac{\left(?\left(\frac{1}{2-x}\right)-\frac{1}{2}\right)\psi_n\left(\frac{1}{2-x}\right)}{(2-x)^2}-\frac{\left(?\left(\frac{1}{1+x}\right)-\frac{1}{2}\right)\psi_n\left(\frac{1}{1+x}\right)}{(1+x)^2}\right)
=\frac{1}{2^n}\left(\frac{\frac{?(x)}{2}\psi_n\left(\frac{1}{2-x}\right)}{(2-x)^2}+\frac{\left(\frac{?(x)}{2}-\frac{1}{2}\right)\psi_n\left(\frac{1}{1+x}\right)}{(1+x)^2}\right)
$$
$$
\begin{aligned}& =\frac{1}{2^{n+1}}\left(?(x)-\frac{1}{2}\right)\psi_{n+1}(x)+\frac{1}{2^{n+2}}\left(\frac{\psi_n\left(\frac{1}{2-x}\right)}{(2-x)^2}-\frac{
\psi_n\left(\frac{1}{1+x}\right)}{(1+x)^2}\right)\\=&\frac{1}{2^{n+1}}\left(?(x)-\frac{1}{2}\right)\psi_{n+1}(x)\\&+\frac{1}{2^{n+2}}\sum_{\frac{p}{q}\in W_n}\left(\frac{1}{((2q-p)-qx)^2}+\frac{1}{((2q-p)-(q-p)x)^2}-\frac{1}{(2q-p)-q(1-x))^2}-\frac{1}{(2q-p)-(q-p)(1-x))^2}\right)\\=&\frac{1}{2^{n+1}}\left(?(x)-\frac{1}{2}\right)\psi_{n+1}(x)+\frac{1}{2^{n+2}}\sum_{\frac{p}{q}\in W_{n+1}}\left(\frac{1}{(q-px)^2}-\frac{1}{(q-p(1-x))^2}\right).\end{aligned}
$$
Therefore, from inductive assumption for $n$  we deduce
$$
h_{n+1}(x)=\frac{1}{2^{n+1}}\left(?(x)-\frac{1}{2}\right)\psi_{n+1}(x)+\frac{1}{2^{n+2}}\sum_{\frac{p}{q}\in W_{n+1}}\left(\frac{1}{(q-px)^2}-\frac{1}{(q-p(1-x))^2}\right)\\
$$
\begin{equation}\label{equa}
+\sum_{\frac{p}{q}\in W_n}\left(c_{\frac{p}{q}}\left(\frac{1}{((2q-p)-qx)^2}-\frac{1}{((q-p)+qx)^2}\right)-
c_{\frac{p}{q}}\left(\frac{1}{((2q-p)-(q-p){x})^2}-\frac{1}{(q+(q-p)x)^2}\right)\right).
\end{equation}.
We should note that  by (\ref{wn+1}) for $\frac{p}{q} \in W_{n+1}$ we have either
\begin{equation}\label{ei1}
\frac{p}{q} 
=
\frac{q'}{2q'-p'} \,\,\,\,\, \text{with}\,\,\,\,\, \frac{p'}{q'} \in W_n
\end{equation}
and
$$
\frac{1}{(q-px)^2}-\frac{1}{(q-p(1-x))^2} =
\frac{1}{((2q'-p')-qx)^2}-\frac{1}{((2q'-p')-q'(1-x))^2}=
\frac{1}{((2q'-p')-qx)^2}-\frac{1}{((q'-p')+q'x)^2}
,
$$
or
\begin{equation}\label{ei2}
\frac{p}{q} 
=
\frac{q'-p'}{2q'-p'}  \,\,\,\,\, \text{with}\,\,\,\,\,  \frac{p'}{q'} \in W_n
\end{equation}
and
$$
\frac{1}{(q-px)^2}-\frac{1}{(q-p(1-x))^2} =
\frac{1}{((2q'-p')-(q'-p')x)^2}-\frac{1}{(q'+(q'-p')x)^2}.$$
So  from equality  (\ref{equa})   we get
$$
h_{n+1}(x)=\frac{1}{2^{n+1}}\left(?(x)-\frac{1}{2}\right)\psi_{n+1}(x)+\sum_{\frac{p}{q}\in W_{n+1}} c_{\frac{p}{q}}\left(\frac{1}{(q-px)^2}-\frac{1}{(q-p(1-x))^2}\right),
$$
where
$$
c_{\frac{p}{q}}=
\begin{cases}
\frac{1}{2^{n+2}}+c_{\frac{p'}{q'}}\,\,\,\text{when (\ref{ei1}) holds},\cr
\frac{1}{2^{n+2}}-c_{\frac{p'}{q'}}\,\,\,\text{when (\ref{ei2}) holds}
\end{cases}
$$
By inductive assumption 
$| c_{\frac{p'}{q'}}|\le\frac{1}{2}-\frac{1}{2^{n+1}} $, and  so
$|c_{\frac{p}{q}}|\le\frac{1}{2}-\frac{1}{2^{n+1}}+ \frac{1}{2^{n+2}} = \frac{1}{2}-\frac{1}{2^{n+2}}$.\\
Hence (\ref{dortmund}) follows by induction.
Now, again applying the differential mean-value theorem to (\ref{dortmund}) proves Lemma \ref{l8}.$\Box$

\vskip+0.3cm

{ Finally, }
  Lemma \ref{l6v} follows from the previous three lemmas.$\Box$

\vskip+0.3cm

   \subsection{Technical lemma}
   
   Consider the quantity
   $$
   \sigma_n(y) =
   \sum_{k=1}^n
   \frac{1}{(y+k)^2} \,
   \min \left(
   \frac{1}{(y+k)\log(n-k+2)}, \frac{1}{(n-k)^2}\right).
   $$
   
   \begin{lem}\label{ltex} 
 $$
 \int_1^\infty \sigma_n(y) {\rm d}y  = O(n^{-\frac{3}{2}}(\log n)^{-\frac{1}{2}}),\,\,\,\,
 n \to \infty.
$$
\end{lem}
   
   {\bf Proof.} We divide  sum $
   \sigma_n(y) $ into three:
   $$
   \sigma_n(y)  = \Sigma^{(1)}+\Sigma^{(2)}+\Sigma^{(3)},\,\,\,\,\,\,
   \Sigma^{(1)}=
   \sum_{k\le \frac{n}{2}},\,\,\,\,\,
   \Sigma^{(2)}
   =
      \sum_{\frac{n}{2}< k\le  n - \delta(n)}
     ,
     \,\,\,\,\,
   \Sigma^{(3)}
   = \sum_{  n - \delta(n)<k\le n},
   $$where $\delta(n)=n^\frac{1}{2}(\log n)^\frac{1}{2}$.\\
   For the first sum we have
   $$
    \Sigma^{(1)}\le 
       \sum_{k\le \frac{n}{2}}
   \frac{1}{(y+k)^2} \cdot
   \frac{1}{(n-k)^2}
   \ll
   \frac{1}{n^2}
   \sum_{k=1}^{n} \frac{1}{(y+k)^2}
   \ll
    \frac{1}{n^2}
    \left( \frac{1}{y} - \frac{1}{y+n}\right),
    $$
    and so
    $$
    \int_1^\infty   \Sigma^{(1)}{\rm d}y\ll
     \frac{1}{n^2}
     \int_1^\infty
    \left( \frac{1}{y} - \frac{1}{y+n}\right){\rm d}y
    \ll
    \frac{\log n}{n^2}.
    $$
    Now we calculate the upper bound for the second sum
    $$
      \Sigma^{(2)} \le
      \sum_{\frac{n}{2}< k\le  n - \delta(n)} 
   \frac{1}{(y+k)^2} \cdot
   \frac{1}{(n-k)^2}
   \ll
   \frac{1}{(y+n)^2}
   \sum_{\delta(n)\le k_1\le n } \frac{1}{k_1^2}\ll
      \frac{1}{(y+n)^2\delta(n)},
    $$
    and so
     $$
    \int_1^\infty   \Sigma^{(2)}{\rm d}y\ll
     \frac{1}{\delta(n)}
     \int_1^\infty
    \frac{{\rm d}y}{(y+n)^2} 
    \ll
    n^{-\frac{3}{2}}(\log n)^{-\frac{1}{2}}.
    $$
    Finally, for the third sum we have
    $$
     \Sigma^{(3)} \le
     {
      \sum_{n-\delta(n)<k\le n}
   \frac{1}{(y+k)^3\log(n-k+2)}
   }
   $$
   $$
        {
   \ll
   \frac{1}{(y+n)^3}
   \sum_{1\le k_1\le \delta (n)}
   \frac{1}{\log(k_1+2)}
\ll
   \frac{\delta(n)}{\log\delta(n)}\cdot\frac{1}{(y+n)^3}
   \ll\frac{n^\frac{1}{2}}{(y+n)^3(\log n)^\frac{1}{2}}  
   },   
   $$
   and 
     $$
    \int_1^\infty   \Sigma^{(3)}{\rm d}y\ll 
    n^\frac{1}{2}(\log n)^{-\frac{1}{2}}     \int_1^\infty
    \frac{{\rm d}y}{(y+n)^3} 
    \ll
    n^{-\frac{3}{2}}(\log n)^{-\frac{1}{2}}.
    $$
    Lemma is proven.$\Box$

\subsection{End of the proof of Theorem \ref{thm1}}
{ Because of (\ref{R_beta}), it would be enough to prove the upper bound
\begin{equation}\label{finn}
\beta_n= \int_{\frac{1}{2}}^1\rho(2^n?(u))f_0(u)\text{d}u  = O(n^{-\frac{3}{2}}(\log n)^{-\frac{1}{2}}),\,\,\,n\to \infty.
\end{equation}
 }
By Lemma \ref{l3} 
{ by means of change of variables
$ u =\frac{y}{y+1}$}
for $\beta_n$ we have
\begin{equation}\label{finn1}
\beta_n=\int_{\frac{1}{2}}^1\rho(?(u))f_n(u)\text{d}u\ll\int_\frac{1}{2}^1|f_n(u)|\text{d}u=\int_1^\infty\frac{1}{(y+1)^2}\left|f_n\left(\frac{y}{y+1}\right)\right|\text{d}y.
\end{equation}
{ Now we combine together all our auxiliary results: 
formula  (\ref{recurrence}) from Lemma \ref{l50},  formula  (\ref{l5_2}) from Lemma \ref{l5}, 
Main Lemma \ref{l6v} and technical Lemma \ref{ltex}.

Namely, (\ref{l50}) gives us}
\begin{equation}\label{finn2}
\int_1^\infty\frac{1}{(y+1)^2}\left|f_n\left(\frac{y}{y+1}\right)\right|\text{d}y=\int_1^\infty\left|
  \frac{f_0
  \left(
  \frac{y+n}{y+n+1}
  \right)}{
  (y+n+1)^2}
   -\sum_{k=1}^n 
   \frac{f_{n-k} 
   \left(
  \frac{y+k}{2y+2k-1}
  \right)}
  {(2y+2k-1)^2}\right|\text{d}y
  .
  \end{equation}
  {
Then as
$?\left(
\frac{y+n}{y+n+1} \right)  = ?\left( 1-
\frac{1}{y+n+1} \right) =1+O(2^{-(y+n)}),$ it is clear that 
\begin{equation}\label{finn3}
f_0\left(
\frac{y+n}{y+n+1} \right) =
?\left( 1-
\frac{1}{y+n+1} \right) - 1 + \frac{1}{y+n+1}= O\left( \frac{1}{y+n}\right),
\end{equation}
meanwhile by  (\ref{l5_2})  
and equalty
$$
\frac{y+k}{2y+2k-1}
=
  \frac{1}{2} -
\frac{1}{2(2y+2k-1)}  
$$
we get
\begin{equation}\label{finn4}
f_{n-k}\left(
\frac{y+k}{2y+2k-1} \right) 
 =
O\left(
\frac{1}{(n-k)^2}
\right).
 \end{equation}
 By Main Lemma  \ref{l6v}  we get
 \begin{equation}\label{finn5}
 f_{n-k}\left(
\frac{y+k}{2y+2k-1} \right)  =
O\left(
\frac{1}{(y+k)\log (n-k+2)}
\right).
 \end{equation}
 Now we substitute (\ref{finn2} - \ref{finn5})  into (\ref{finn1}) and  by Lemma \ref{ltex} obtain 
  $$
  \beta_n
  \ll\int_1^\infty\left(\frac{1}{(y+n)^3}+\sigma_n(y)\right)\text{d}y=O(n^{-\frac{3}{2}}(\log n)^{-\frac{1}{2}}).$$  Inequality (\ref{finn}) and so} Theorem \ref{thm1} is proven.$\Box$

  \section{Proof of Theorem \ref{thm2}}\label{vvv}

  By  Lemma 1 from \cite{g}, Theorem   \ref{thm2} is equivalent to 
  $$\alpha_n:=\int_{\frac{1}{2}}^1\cos(2\pi?(x))\psi_n(x)\text{d}x=\frac{1}{2}+O\left(\frac{1}{\log n}\right).$$
  By the explicit formula (\ref{eqn-1})  and the monotonicity of $?(x)$ we have $$1-?(x)=?(1-x)\le?\left(\frac{1}{\left\lfloor\frac{1}{1-x}\right\rfloor}\right)=2^{-\left\lfloor\frac{1}{1-x}\right\rfloor+1}\ll2^{-\frac{1}{1-x}}$$
  uniformly in 
 $x\in (0,1)$.
  
 From the explicit formula (\ref{eqn-30}) and the result  by
 Kesseböhmer  and Stratmann (\ref{kss}) from \cite{ks}   we deduce that  
 $$\psi_n(x)=\sum_{\frac{p}{q}\in W_n}\left(\frac{1}{(q-px)^2}+\frac{1}{(q-p(1-x))^2}\right)\le S_n\cdot\left(\frac{1}{(1-x)^2}+\frac{1}{x^2}\right)\ll\frac{1}{(1-x)^2\log n}$$
   uniformly in 
 $x\in \left[\frac{1}{2},1\right)$.
 Now
 $$\frac{1}{2}-\alpha_n=\int_{\frac{1}{2}}^1\left(1-\cos(2\pi?(x))\right)\psi_n(x)\text{d}x\ll\int_{\frac{1}{2}}^1(1-?(x))^2\psi_n(x)\text{d}x\ll\frac{1}{\log n}\int_{\frac{1}{2}}^12^{-\frac{2}{1-x}}\cdot\frac{1}{(1-x)^2}\text{d}x \ll\frac{1}{\log n}.$$
 Everything is proven.$\Box$

 \vskip+0.3cm
 
\begin{large}\begin{bf}Acknowledgements\end{bf}\end{large}\,\,
The authors would like to thank Professor \begin{bf}Nikolay Moshchevitin\end{bf} for suggesting this project, giving us inspiring hints, providing helpful references and helping make our original version more readable. We also thank Professor \begin{bf}Oleg N. German\end{bf} as well as \begin{bf}Jianqiao Xia\end{bf} and \begin{bf}Marsault Chabat\end{bf} for checking our proofs, giving us suggestions and helping us to find references.

\end{document}